\documentclass[11pt]{amsart}
\usepackage{amsmath,amsthm}
\usepackage{amssymb,latexsym}
\usepackage{enumerate}

\newtheorem{thm}{Theorem} %[section]

\theoremstyle{definition}
\newtheorem{rem}{Remark}

\numberwithin{equation}{section}
\makeatletter
\@namedef{subjclassname@2010}{%
  \textup{2010} Mathematics Subject Classification}
\makeatother

\newcommand{\Q}{\mathbb{Q}}
\newcommand{\Z}{\mathbb{Z}}

\newcommand{\N}{\mathbb{N}}
\newcommand{\Pb}{\mathbb{P}}

\newcommand{\Oc}{\mathcal{O}}

\begin{document}

%%%%%%%%%%%%%%%%

\title[On $k$-free values]{On the $k$-free values of the polynomial $xy^k+C$}

\author[K. Lapkova]{Kostadinka Lapkova}
\address{MTA Alfr\'ed R\'enyi Institute of Mathematics\\
1053 Budapest, Re\'altanoda u. 13-15, Hungary}
\email{lapkova.kostadinka@renyi.mta.hu}

%\author[K. Lapkova]{Kostadinka Lapkova\\
%MTA Alfr\'ed R\'enyi Institute of Mathematics\\
%1053 Budapest, Re\'altanoda u. 13-15, Hungary}\\
%Email:\texttt{lapkova.kostadinka@renyi.mta.hu}

\date{20.10.2015}

\begin{abstract} Consider the polynomial $f(x,y)=xy^k+C$ for $k\geq 2$ and any nonzero integer constant $C$. We derive an asymptotic formula for the $k$-free values of $f(x,y)$ when $x, y\leq H$. We also prove a similar result for the $k$-free values of $f(p,q)$ when $p,q\leq H$ are primes, thus extending Erd\H{o}s' conjecture for our specific polynomial. The strongest tool we use is a recent generalization of the determinant method due to Reuss.
\end{abstract}

% 11N32:other multiplicative structure of polynomial values; 11N37 : asymptotic results on arithmetic functions 


\maketitle

\section{Introduction}\label{sec:intro}

\hspace{0.5cm}

Let $k$ and $n$ be integers and $k\geq 2$. We say that $n$ is \emph{$k$-free} if there is no prime $p$ such that $p^k\mid n$. Consider the irreducible polynomial $f(n)\in\Z\left[x\right]$ of degree $d$. Let also $f(n)$ have no fixed \emph{$k$-th power prime divisor}, i.e. there is no prime $p$ for which $p^k\mid f(n)$ for every $n\in\Z$. Then we expect that the set $f(\Z)=\left\{f(n), n\in\Z\right\}$ contains infinitely many $k$-free values. The first one who obtained a result in this direction was Ricci \cite{ricci}, who derived an asymptotic formula for the quantity
\begin{equation}\label{quant1} \#\left\{n\leq H : f(n) \text{ is }k\text{-free}\right\}
\end{equation}
when $k\geq d$. Later, Erd\H{o}s \cite{erdos} proved the conjecture in the case $k=d-1$ and $d\geq 3$, and Hooley \cite{hooley67} obtained an asymptotic formula for each such $k$. Some of the recent results on this conjecture are asymptotic formulas for the $k$-free values due to Browning \cite{brown} when $k\geq (3d+1)/4$ and $d\geq 3$, and to Heath-Brown \cite{heath-br2} for the polynomial $f(x)=x^d+C$ when $k\geq(5d+3)/9$ for $d\geq 3$. In general the problem is harder with $d$ increasing, but easier the closer $k$ is to $d$.\\

Similarly for irreducible multi-variable polynomials $f(t_1,\ldots,t_r)$ of degree $d$, which do not have a fixed $k$-th power prime divisor, we conjecture that the quantity 
\begin{equation}\label{quantK} \#\left\{t_1,\ldots,t_r\leq H : f(t_1,\ldots,t_r) \text{ is }k\text{-free}\right\}
\end{equation}
has a positive density. Recent results in this direction are due to Tolev \cite{tolev} and Le Boudec \cite{leBoudec} who investigate carefully asymptotic formulas for specific polynomials. Earlier result of Poonen \cite{poonen}, which lays on the $abc$ conjecture, shows that the number of the square-free values of any $f(t_1,\ldots,t_r)$ has a positive density, which, however, is differently defined than the traditional density which arises when we evaluate asymptotically (\ref{quantK}).\\

Naturally the most extensively studied multi-variable polynomials are the two-variable polynomials, where the situation differs in the homogeneous and inhomogeneous case. While the power-free values of binary forms were investigated already by Greaves \cite{greaves}, it is only recently that results were obtained for inhomogeneous two-variable polynomials by Hooley \cite{hooley09}. Later Browning \cite{brown} achieved the smallest $k>39d/64$ up to date. Note that the latter papers derive not an asymptotic formula but a positive lower bound of the expected true order of magnitude for the quantity (\ref{quantK}). Nevertheless, in the case $k=d-1$ for two-variable polynomials, as in Hooley's papers \cite{hooley67} and \cite{hooley76} for one-variable polynomials, it is expected that we can derive an asymptotic formula.\\

In this paper we consider the inhomogeneous polynomial $f(x,y)=xy^k+C$ for $k\geq 2$ and any nonzero integer constant $C$. We derive an asymptotic formula for its $k$-free values when $x,y\leq H$.

\begin{thm}\label{thm:L1}Let $f(x,y)=xy^k+C\in\Z[x,y]$ for $k\geq 2$ and $C\neq 0$. Let $S(H)$ count the $k$-free values of $f(x,y)$ when $1\leq x,y\leq H$. Then, for some real $\delta=\delta(k,f)>0$, we have
$$S(H)=c_{f}H^2+\Oc\left(H^{2-\delta}\right),$$
where 
$$c_{f}=\prod_p \left(1-\frac{\rho(p^k)}{p^{2k}}\right)$$
and 
$$\rho(m)=\# \left\{ (\mu,\nu)\in(\Z/m\Z)^2 : \quad f(\mu,\nu) \equiv 0\pmod m \right\}.$$
\end{thm}

\begin{rem} We will show that the infinite product defining $c_{f}$ is convergent and is not zero. For $k\geq 3$ we can take $\delta=1/(7k)$ and for $k=2$ we can take $\delta=2-\varepsilon-G_2$, where $G_k$ is the expression (\ref{eq:G}) and $\varepsilon>0$ is arbitrary small. This way for $k=2$ we can get an error term $\Oc\left(H^{1.979}\right)$. We do not guarantee that this way $\delta$ is the best possible, but it is close to optimal by the choices made during the proof. 
\end{rem}

The determinant method of Heath-Brown \cite{heath-br0} has been already successfully used for counting power-free values of polynomials in one variable by Heath-Brown himself in  \cite{heath-br2} and \cite{heath-br1}, by Browning \cite{brown} and Reuss \cite{reuss2}. In this paper we apply a major result of Reuss \cite{reuss1} which generalizes the approximate determinant method developed by Heath-Brown in \cite{heath-br1}. As we reduce our problem of counting $k$-free values of the polynomial $f(x,y)=xy^k+C$ to the problem of counting solutions of a Diophantine equation inside a certain box, we can directly use Theorem 1 of Reuss \cite{reuss1}, supplying one more application of his powerful result. Note that for our particular polynomial it would suffice to use Lemma $5$ from Dietmann-Marmon \cite{diet-mar}, where again by determinant method they give upper bound of the solutions of the equation $ax^k-by^k=1$ with restricted sizes of the variables. Reuss however resolves the more general equation $a^lx^k-b^ly^k=C$ for any integers $k,l,C$ with $1\leq l<k$ and $C\neq 0$.\\ 

In \cite{erdos} Erd\H{o}s conjectured that similarly the set $f(\Pb)=\left\{f(p), p  \text{ prime}\right\}$ for a polynomial of degree $d$ with no fixed $(d-1)$-th power prime divisor contains infinitely many $(d-1)$-free values. Heath-Brown's strategy for finding an asymptotic formula for the quantity (\ref{quant1}) was successfully used for this problem as well, e.g. in \cite{brown}, \cite{heath-br2}. Finally Reuss \cite{reuss2} settled the conjecture for all $d\geq 3$, demonstrating the strength of this (mostly) analytic method over other approaches, e.g. of Helfgott \cite{helfgott}.\\

We extend Erd\H{o}s' conjecture for the set $f(\Pb, \Pb)=\left\{f(p,q),\, p\text{ and }q \text{ primes}\right\}$ for our specific two-variable polynomial. 

\begin{thm}\label{thm:L2}Let $f(x,y)=xy^k+C\in\Z[x,y]$ for $k\geq 2$ and $C\neq 0$. Let $S'(H)$ count the $k$-free values of $f(p,q)$ for prime numbers  $1< p,q\leq H$. Then, for any real $K>2$, we have the asymptotic formula
$$S'(H)=c'_{f}\pi(H)^2+\Oc\left(\frac {H^2}{(\log H)^K}\right),$$
where 
$$c'_{f}=\prod_p \left(1-\frac{\rho '(p^k)}{\varphi(p^k)^2}\right)$$
and 
$$\rho'(m)=\# \left\{ (\mu,\nu)\in(\Z/m\Z)^2 : \quad (\mu,m)=(\nu,m)=1 \text{ and } f(\mu,\nu) \equiv 0\pmod m \right\}.$$
\end{thm}
\begin{rem}We will show that the constant $c'_{f}$ makes sense, i.e. the infinite product is convergent, and $c'_{f}\neq 0$.
\end{rem}
The motivation for considering power-free values of the polynomial $f(x,y)=xy^k+C$ comes from our previous works \cite{biroL} and \cite{lapkova}. There we resolve the class number one problem for real quadratic fields with discriminant $(an)^2+4a$. It is natural to ask if the polynomial $(an)^2+4a$ takes square-free values with positive density. We answer affirmatively this question for the irreducible factor $f(x,y)=xy^2+4$. Actually with Theorem \ref{thm:L2} we give an asymptotic formula for the square-free values of $(xy)^2+4x$ when both $x$ and $y$ are primes, presumably a harder problem than the one when both factors $x$ and $xy^2+4$ are square-free.

\section{Notations}
\hspace{0.5cm}

We reserve the letters $p$ and $q$ for primes and $\varepsilon$, $\epsilon$ for arbitrary small positive numbers, not necessarily the same in different occasions. By $\left[t\right]$ we denote the integer part of the real number $t$, and $(a,b)$ denotes the greatest common divisor of the integers $a$ and $b$. By $\bar y$ we denote the multiplicative inverse of $y$ modulo some integer which should be clear by the context. The dependences in the Landau symbol $\Oc$ and the Vinogradov symbol $\ll$ usually will be omitted.  Sometimes in the summation sign we write only the variables which we sum, but again the intervals they lie in should be clear by the context.Sometimes we write $(m)$ instead of $\pmod m$. %We use the usual $e(x)=e^{2\pi }$
\par By $\pi(H)$ as usual we denote the prime-counting function, and by $\pi(H;q,a)$ the prime-counting function in the arithmetic progression $\{a+nq, n\geq 0\}$. The Euler function is denoted by $\varphi(n)$. By $\tau(n)$ we denote the number of divisors function. We write $x\sim X$ to say that $X< x< 2X$ and we write $x\asymp X$ to say that there exist positive constants $A,B$, independent of $X$, such that $AX\leq |x|\leq BX$.
%structure of the paper
\par In the next two sections we prove Theorem \ref{thm:L1}, where the methods used in \S\ref{sec:3} are elementary and in \S\ref{sec:4} we formulate and apply Reuss' theorem. We prove Theorem \ref{thm:L2} in section \S\ref{sec:5} .

\section{Proof of Theorem \ref{thm:L1}: elementary tools}\label{sec:3}

\hspace{0.5cm} Let $$\mu_k(n)=\left\{
\begin{array}{lll}1 & , & \text{if n is k-free}\\
0 & , & \text{otherwise}
\end{array}\right.$$ 
be the characteristic function of the $k$-free numbers. Just like with the M\"obius function we have the identity
\begin{equation}\label{eq:mu} \mu_k(n)=\sum_{d^k\mid n}\mu(d)\,.
\end{equation}
In order to count the $k$-free values of the polynomial $f(x,y)=xy^k+C$, where $k\geq 2$ is an integer and $C$ is a non-zero integer, and $1\leq x,y \leq H$, we can consider the sum
\begin{equation}S:=S(H)=\sum_{1\leq x,y \leq H}\mu_k(xy^k+C)\,.
\end{equation}
Using the equation (\ref{eq:mu}) and changing the order of summation we transform the sum $S$.
$$S=\sum_{1\leq x,y \leq H}\sum_{d^k\mid f(x,y)}\mu(d)=\sum_{1\leq d\leq f(H,H)^{1/k}} \mu(d) \sum_{\substack{1\leq x,y\leq H \\ f(x,y)\equiv 0\pmod{d^k}}}1=\sum_{1\leq d\ll H^{(k+1)/k}}\mu(d)S(d^k,H),$$
where we denote 
\begin{equation}S(m,H):=\sum_{\substack{1\leq x,y\leq H\\ f(x,y)\equiv 0\pmod m}}1\,.
\end{equation}

We split the sum $S$ into four parts: 
\begin{equation}\label{eq:S_split}S=S_1+S_2+S_3+S_4,
\end{equation}
where $$S_1:=\sum_{1\leq d\leq z_1}\mu(d)S(d^k,H)$$
for a parameter $z_1$, and the summation is correspondingly over $z_1<d\leq z_2$ for $S_2$, $z_2<d\leq z_3$ for $S_3$ and $z_3<d\ll H^{(k+1)/k}$ for $S_4$. The parameters $z_1, z_2, z_3$ are to be chosen soon in a convenient way. Essentially the estimates of $S_1, S_2$ and $S_4$ will be trivial, and for $S_3$ we will apply Theorem $1$(Reuss, \cite{reuss1}).\\

%%----------------------------------------------------------------------

\subsection{Estimation of $S_1$} The sum $S_1$ will contribute the main term in our asymptotic formula. Introduce the notation
\begin{equation} \label{eq:M} M(\alpha, m, H):=\sum_{\substack{1\leq x\leq H\\ x\equiv\alpha\pmod m}} 1\,.
\end{equation}
Then 
$$S(m,H)=\sum_{\substack{1\leq \mu,\nu \leq m \\ f(\mu,\nu)\equiv 0(m)}}M(\mu,m,H)M(\nu,m,H)\,.$$
Clearly 
%
%$$M(\alpha,m,H)=\sum_{\substack{x\leq H \\ x\equiv \alpha (m)}}1=\left[\frac H m\right]+\Oc(1)=\frac H m +\Oc(1)\,,$$
%
$$M(\alpha,m,H)=\frac H m +\Oc(1)\,,$$
therefore 
%\begin{eqnarray}\label{eq:SmH} S(m,H)&=&\sum_{\substack{\mu,\nu\leq m\\ f(\mu,\nu)\equiv 0(m)}}\left(\frac H m +\Oc(1)\right)^2=\sum_{\substack{\mu,\nu\leq m\\ f(\mu,\nu)\equiv 0(m)}}\left(\frac{H^2}{m^2} +\Oc\left(\frac H m\right)+\Oc(1)\right)\nonumber\\
%&=& H^2\frac{\rho(m)}{m^2}+\Oc\left(H\frac{\rho(m)}{m}\right)+\Oc(\rho(m))\,,
%\end{eqnarray}
\begin{equation}\label{eq:SmH} S(m,H) = H^2\frac{\rho(m)}{m^2}+\Oc\left(H\frac{\rho(m)}{m}\right)+\Oc(\rho(m))\,,
\end{equation}
where $\rho(m)$ was defined in Theorem \ref{thm:L1}.
Then 
\begin{eqnarray}\label{eq:S1}S_1 & = & \sum_{1\leq d\leq z_1}\mu(d)\left[H^2\frac{\rho(d^k)}{d^{2k}}+\Oc\left(H\frac{\rho(d^k)}{d^k}\right)+\Oc\left(\rho(d^k)\right)\right]\nonumber\\
& =  & H^2\sum_{1\leq d\leq z_1}\mu(d)\frac{\rho(d^k)}{d^{2k}}+\Oc\left(H\sum_{d\leq z_1}\frac{\rho(d^k)}{d^k}\right)+\Oc\left(\sum_{d\leq z_1}\rho(d^k)\right)\nonumber\\
& = & H^2\sum_{d=1}^\infty\mu(d)\frac{\rho(d^k)}{d^{2k}} -H^2\sum_{d > z_1}\mu(d)\frac{\rho(d^k)}{d^{2k}} +\Oc\left(H\sum_{d\leq z_1}\frac{\rho(d^k)}{d^k}\right)+\Oc\left(\sum_{d\leq z_1}\rho(d^k)\right)\,.
\end{eqnarray}
Hooley's Lemma $1$ from \cite{hooley09} claims that the function $\rho(m)$ is multiplicative and $\rho(p^k)=\Oc(p^{2k-2})$ for $k\geq 2$. Therefore for square-free $d$ we have $\rho(d^k)=\Oc(d^{2k-2+\epsilon})$ for $\epsilon>0$ arbitrary small. Then the infinite sum
\begin{equation}c_{f}:=\sum_{d=1}^\infty\mu(d)\frac{\rho(d^k)}{d^{2k}}=\prod_p \left(1-\frac{\rho(p^k)}{p^{2k}}\right)
\end{equation}
is convergent because its common term $\mu(d)\rho(d^k)/d^{2k}=\Oc(d^{2k-2+\epsilon})/d^{2k}=\Oc(d^{-2+\epsilon})$. Something more, $c_{f}\neq 0$. Indeed, $c_{f}=0$ only if it has a zero factor, i.e. $\rho(p^k)=p^{2k}$ for some $p$. Take some $\sigma\in \Z/p^k\Z$ such that $\sigma\not\equiv -C\pmod{p^k}$. Then $f(\sigma,1)=\sigma\cdot 1^k+C\not\equiv 0\pmod{p^k}$, so surely $\rho(p^k)<p^{2k}$.\\

Further, using that in all summands in (\ref{eq:S1}) the variable $d$ is square-free, thus again $\rho(d^k)=\Oc(d^{2k-2+\epsilon})$, we get
$$\sum_{d > z_1}\mu(d)\frac{\rho(d^k)}{d^{2k}}\ll\sum_{d>z_1}d^{-2+\epsilon}\ll z_1^{-1+\epsilon}\,,$$
$$\sum_{d\leq z_1}\frac{\rho(d^k)}{d^k}\ll\sum_{d\leq z_1}d^{k-2+\epsilon}\ll z_1^{k-1+\epsilon}\,,$$
\begin{equation}\label{eq:ro}\sum_{d\leq z_1}\rho(d^k)\ll\sum_{d\leq z_1}d^{2k-2+\epsilon}\ll z_1^{2k-1+\epsilon}\,.
\end{equation}
Putting all these into (\ref{eq:S1}) we get
\begin{equation}\label{eq:S1err}S_1=c_{f}H^2+\Oc(H^2z_1^{-1+\epsilon})+\Oc(Hz_1^{k-1+\epsilon})+\Oc(z_1^{2k-1+\epsilon})\,.
\end{equation}

Let us choose an appropriate value for $z_1$. Write $z_1=H^\alpha$ for some real number $0<\alpha<1$. For the first error term in (\ref{eq:S1err}) we always have $2-\alpha+\epsilon<2$ as long as $\epsilon>0$ is small enough. In the second error term we want to have $1+\alpha(k-1+\epsilon)<2$,  i.e. $\alpha<1/(k-1+\epsilon)<1/(k-1)$. In the third error term we require $\alpha(2k-1+\epsilon)<2$, i.e. $\alpha<2/(2k-1+\epsilon)<1/(k-1)$. Then $\alpha=1/k$ is a good choice and we get
$$S_1=c_{f}H^2+\Oc(H^{2-1/k+\epsilon})+\Oc(H^{1+(k-1)/k+\epsilon})+\Oc(H^{(2k-1)/k+\epsilon})\,,$$
or finally
\vspace{0.2cm} 
\begin{equation}\label{eq:S1final}   
S_1=\sum_{1\leq d\leq H^{1/k}}\mu(d)S(d^k,H)=c_{f}H^2+\Oc\left(H^{2-1/k+\epsilon}\right)\,.
\end{equation}

%%----------------------------------------------------------------------

\subsection{Estimation of $S_2$}\label{sec:S2} It turns out that for the estimation of the second sum trivial arguments yield better results than more elaborate tools like exponential sums and their Weil-type of estimates, e.g. like the ones used in \cite{leBoudec} and \cite{tolev}. The gain comes if we straight choose the parameter $z_2=H^{1-\delta}$ for a small $\delta>0$, to be specified later. Then 
$$S_2=\sum_{H^{1/k}<d\leq H^{1-\delta}}\mu(d)\sum_{\substack{x,y\leq H\\ xy^k+C\equiv 0(d^k)}}1\ll \sum_{H^{1/k}<d\leq H^{1-\delta}}\sum_{\substack{x,y\leq H\\ xy^k\equiv -C(d^k)}}1\,.$$

%Suppose first that $(y,d)=1$. If we fix $d$ and $y$ the innermost sum counts those $x\leq H<d^k$ which satisfy the congruence $x\equiv -C\bar{y}^k\pmod {d^k}$. But $x$ runs values in only one copy of a full residue system modulo $d^k$, so there is at most one solution of the latter congruence. Using the notation (\ref{eq:M}) we have 
%$$ S_2^{*}:= \sum_{H^{1/k}<d\leq H^{1-\delta}}  \sum_{\substack{y\leq H \\ (y,d)=1 }} M(-C\bar{y}^k,d^k,H)\ll H^{1-\delta}\cdot H\ll H^{2-\delta}\,.$$ 
 If we fix $d$ and $y$, the innermost sum counts those $x\leq H<d^k$ which satisfy the congruence $xy^k\equiv -C\pmod {d^k}$. This congruence is solvable in $x$ if and only if $(y^k,d^k)\mid C$, and in this case in a full residue system modulo $d^k$ there are exactly $(y^k,d^k)$ solutions. If $(y,d)^k\mid C$, we have $(y,d)^k\ll 1$. Then 
\begin{equation}\label{eq:S2final}S_2\ll \sum_{H^{1/k}<d\leq H^{1-\delta}}  \sum_{\substack{y\leq H \\ (y,d)^k\mid C }} (y,d)^k\ll H^{1-\delta}\cdot H\ll H^{2-\delta}\,.
\end{equation}

%If $(y,d)=\Delta>1$, then it is necessary that $\Delta^k$ divides $C$. Clearly if $\mu_k(C)=1$ this does not happen at all, and in case $C$ has a $k$-th power divisor, there are not more than $\tau(C)=\Oc(1)$ of these divisors. Then the solutions of the congruence  
%\begin{equation}\label{eq:congr}xy^k\equiv -C\pmod{d^k}
%\end{equation} 
%for fixed $d$ and $y$ are $\Delta^k=\Oc(1)$, and the upper bound for $S_2^{*}$ applies in this case too. Then we have
%\begin{equation}\label{eq:S2final}
%S_2=\sum_{H^{1/k}\leq d\leq H^{1-\delta}}\mu(d)S(d^k,H)=\Oc(H^{2-\delta})\,.
%\end{equation}

%%----------------------------------------------------------------------

\subsection{Estimation of $S_4$} To estimate the tail of the sum $S$ we will use a well-known upper bound for the number of divisors function, i.e. for any positive $\varepsilon>0$ we have $\tau(n)\ll n^\varepsilon$ (Theorem $315$, \cite{hardy-wright}). Then
\begin{eqnarray*} S(d^k, H) & = & \sum_{\substack{x,y\leq H \\ xy^k\equiv -C(d^k)}}1\leq \sum_{\substack{n\ll H^{k+1}\\ n\equiv -C(d^k)}}\sum_{n=xy^k}1 \ll  \sum_{\substack{n\ll H^{k+1}\\ n\equiv -C(d^k)}}\tau(n)\\
 & \ll & H^{(k+1)\varepsilon}  \sum_{\substack{n\ll H^{k+1}\\ n\equiv -C(d^k)}}1\ll H^{\varepsilon}\frac{H^{k+1}}{d^k}=\frac 1 {d^k}H^{k+1+\varepsilon}\,.
\end{eqnarray*}
Then for $S_4$ we get
\begin{eqnarray*} S_4  = \sum_{z_{3} <d\ll H^{1+1/k}}\mu(d)S(d^k, H) & \ll & \sum_{z_3 <d\ll H^{1+1/k}} S(d^k, H) \ll \sum_{z_3<d\ll H^{1+1/k}}\frac 1 {d^k}H^{k+1+\varepsilon}\\
& \ll & H^{k+1+\varepsilon}\sum_{z_3<d}\frac 1 {d^k}\ll H^{k+1+\varepsilon}z_3^{-k+1}\,.
\end{eqnarray*}
%Let $z_3=H^\gamma$ for some real $1\leq \gamma< 1+1/k$. We want to have $H^{k+1+\varepsilon+\gamma(1-k)}<H^2$, i.e. $k-1+\varepsilon<\gamma(k-1)$ or $\gamma > 1+\varepsilon/(k-1)$.\\ 
%
%Set temporarily $z_3=H^{1+\varepsilon^{*}}$, for some $0<\varepsilon^{*}<1/k$. Then we need $k+1+\varepsilon+(1+\varepsilon^{*})(1-k)<2$, i.e. $k+1+\varepsilon-k+1+\varepsilon^{*}(1-k)=2+\varepsilon-\varepsilon^{*}(k-1)<2$ so it is necessary that $\varepsilon^{*}>\varepsilon/(k-1)$, where $\varepsilon$ was arbitrary small and coming from the bound of the divisor function. If we choose $\varepsilon^{*}=2\varepsilon/(k-1)$, the error term turns into $H^{2-\varepsilon}$ with the condition $\varepsilon<(k-1)/(2k)$ which can be always satisfied for small enough $\varepsilon>0$.\\
%
%Finalizing the notation, let assume that 

We set $z_3=H^{1+\delta}$, and for simplicity take this $\delta<1/k$ the same as in the definition of $z_2=H^{1-\delta}$, also choose $\varepsilon=\delta(k-1)/2$. Then we have the upper bound
\begin{equation}\label{eq:S4final}
%S_4=\sum_{H^{1-\delta}<d\ll H^{1+1/k}}\mu(d)S(d^k,H)\ll H^{2-(k-1)\delta/2}\,.
S_4\ll H^{2-(k-1)\delta/2}\,.
\end{equation}

%%----------------------------------------------------------------------

\subsection{Estimation of $S_3$} We will narrow the intervals of summation in $S_3$ by extracting more trivial estimates. This would make possible the application of Theorem $1$(Reuss, \cite{reuss1}), which still does not apply directly to the sum 
\begin{equation} S_3=\sum_{H^{1-\delta}<d\leq H^{1+\delta}}\mu(d)S(d^k, H)=\sum_{H^{1-\delta}<d\leq H^{1+\delta}}\mu(d)\sum_{\substack{1\leq x,y \leq H\\ xy^k+C\equiv 0(d^k)}}1\,.
\end{equation}
Introduce the notation
\begin{equation}
S_3(K_1,K_2;L_1,L_2):=\sum_{H^{1-\delta}<d\leq H^{1+\delta}}\sum_{K_1\leq x\leq K_2}\sum_{\substack{L_1\leq y\leq L_2\\xy^k+C\equiv 0(d^k)}}1\,.
\end{equation}
Clearly $S_3\ll S_3(1,H;1,H)$. 

\subsubsection{The sum $S_3(1,H;1,H^{1-\theta})$} Let us first estimate the contribution to $S_3$ when the variable $y$ lies in the interval $\left[1, H^{1-\theta}\right)$ for some $\theta>0$. We fix $d$ and $y$ and want to solve the congruence $xy^k\equiv -C\pmod{d^k}$, in $x$, for $1\leq x\leq H$. Note that $d^k>H^{k-k\delta}>H$ if $\delta<1-1/k$. The latter holds as we already want $\delta<1/k$ and for $k\geq 2$ we have $1/k\leq 1-1/k$. So, just like in \S\ref{sec:S2}, we will have $\Oc\left((y,d)^k\right)=\Oc(1)$ solutions in $x$. Then 
\begin{equation}\label{eq:S31}
S_3(1,H;1,H^{1-\theta})\ll \sum_{H^{1-\delta}<d\leq H^{1+\delta}}\sum_{1\leq y<H^{1-\theta}}1\ll H^{2-(\theta-\delta)},
\end{equation}
where clearly we must choose $1>\theta>\delta$\,.

\subsubsection{The sum $S_3(1,H^\eta;H^{1-\theta},H)$} Let $0<\eta<1$ be a parameter which we will choose later. We are looking closely at the contribution to the sum $S_3$ when $1\leq x\leq H^\eta$.  We can write the congruence $xy^k+C\equiv 0\pmod{d^k}$ in the equation form 
\begin{equation} \label{eq:equat}xy^k+C=ad^k
\end{equation} for a positive integer $a$. As $ad^k=xy^k+C\ll H^\eta H^k$ and $d^k>H^{k(1-\delta)}$, we should have $a\ll H^{\eta+k}/H^{k(1-\delta)}=H^{\eta+k\delta}$. \\

Let us now fix $x$ and $d$. Observe that the solutions in $a$ of the congruence $ad^k\equiv C\pmod x$ are more than the solutions of the original equation (\ref{eq:equat}). Also the number of solutions is $(d^k,x)\ll 1$ in case $(d^k,x)\mid C$. Then we can bound from above in the following way.
\begin{eqnarray}\label{eq:S3dx}S_3(1,H^\eta;H^{1-\theta},H) & = & \sum_{H^{1-\delta}<d\leq H^{1+\delta}}\sum_{1\leq x\leq H^\eta}\sum_{\substack{y,a\\ ad^k=xy^k+C}}1\ll \sum_{d,x}\sum_{\substack{1\leq a\ll H^{\eta+k\delta}\\ ad^k\equiv C(x)}}1 \nonumber\\
& \ll &  \sum_{d\leq H^{1+\delta}}\sum_{x\leq H^\eta}\sum_{\substack{a\ll H^{\eta+k\delta}\\(d^k,x)\mid C}}(d^k,x)\ll \sum_{d\leq H^{1+\delta}}\sum_{x\leq H^\eta}\sum_{a\ll H^{\eta+k\delta}}1 \,. 
\end{eqnarray}

We get
\begin{equation}\label{eq:S32}
S_3(1,H^\eta;H^{1-\theta},H)\ll H^{1+2\eta+(k+1)\delta} ,
\end{equation}
as long as $\eta>0$ is chosen in such a way that $2\eta+(k+1)\delta<1$. From this condition we also need to have $\delta<1/(k+1)$.

%%----------------------------------------------------------------------

\section{Finishing the proof of Theorem \ref{thm:L1}  : application of Reuss' theorem}\label{sec:4}

In this section we apply Theorem $1$ of Reuss \cite{reuss1} which generalizes the approximate determinant method developed by Heath-Brown in \cite{heath-br1}. % Note that for our purposes we could have used also Lemma 5 from Dietmann-Marmon's work \cite{diet-mar}, but Reuss' result is much more general and our theorem is only one of its many applications, requiring estimation of the solutions of a certain Diophantine equation $ v^le^k-u^ld^k=h$.
 For sake of completeness, and for sake of scrutinizing the conditions of the theorem, we formulate it here.

\begin{thm}[Reuss, 2014]\label{thm:reuss} Let $D,E,z>1$ and $\varepsilon>0$. Let $k,l,h$ be integers such that $1\leq l<k$ and $h\neq 0$. Let
$$\mathcal{N}(z;D,E):=\#\left\{(d,e,u,v)\in\N^4 : d\sim D, e\sim E, u\sim U, v\sim V, v^le^k-u^ld^k=h\right\},$$
where 
$$U=\frac{z^{1/l}}{D^{k/l}}, \text{ and }V=\frac{z^{1/l}}{E^{k/l}}\,.$$
Let $M>1$ be defined by
$$\log M=\frac{9}{8}\frac{\log(DE)\log(UV)}{\log z}\,,$$
and suppose the following conditions are satisfied:
\begin{enumerate}
\item $\log(DE)\asymp \log(UV)\asymp\log z$;
\item $l\geq 2$, or $DE\gg_{k,l,h} z^{1/k}$.
\end{enumerate}
Then, if $z$ is large enough in terms of $\varepsilon$,
$$\mathcal{N}(z;D,E)\ll_{\varepsilon,k,l,h} z^\varepsilon\min\left\{(DEM)^{1/2}+D+E, (UVM)^{1/2}+U+V\right\}\,.$$

\end{thm}

We write $\widetilde{S}_3:=S_3(H^\eta,H;H^{1-\theta},H)$. Putting together the estimates (\ref{eq:S1final}), (\ref{eq:S2final}), (\ref{eq:S4final}), (\ref{eq:S31}), (\ref{eq:S32}) we get
\begin{eqnarray}\label{eq:SOBL} S(H)=c_{f}H^2+\Oc(\widetilde{S}_3)& + & \Oc(H^{2-1/k+\epsilon})+\Oc(H^{2- \delta})+\Oc(H^{2-(k-1)\delta/2})\nonumber\\ & +  &\Oc(H^{2-(\theta-\delta)})+\Oc(H^{1+2\eta+(k+1)\delta})\,,
\end{eqnarray}
and now we examine the last sum $\widetilde{S}_3$.\\

Note that by dyadic subdivision we can write
\begin{eqnarray}\label{eq:dyadic} \widetilde{S}_3 & \ll  & \sum_{H^{1-\delta}<d\leq H^{1+\delta}}\sum_{H^\eta<x\leq H}\sum_{H^{1-\theta}<y\leq H}\sum_{xy^k+C=ad^k}1\nonumber\\
& \ll & (\log H)^3\max_{D,X,Y,A}\sum_{\substack{d\sim D\\x\sim X}}\sum_{\substack{y\sim Y\\a\sim A}}\sum_{xy^k-ad^k=-C}1\nonumber\\
& =&  (\log H)^3\max_{D,X,Y,A} \mathcal{N}(Z;D,Y)\,,
\end{eqnarray}
where $Z=XY^k=AD^k$ and we have 
\begin{eqnarray}\label{eq:Det_int}H^{1-\delta}< & D & \leq H^{1+\delta},\nonumber\\
H^\eta< & X & \leq H,\nonumber\\
H^{1-\theta}< & Y & \leq H,\nonumber\\
1\leq &A&\leq H^{k+1}/D^k\leq H^{1+k\delta}\,.
\end{eqnarray}

We assure that the conditions of Theorem \ref{thm:reuss} hold. First, we should see whether the first condition $\log(DY)\asymp\log(AX)\asymp\log Z$ holds. Using (\ref{eq:Det_int}) we consecutively obtain $H^{\eta+k(1-\theta)}\leq XY^k=Z\leq H^{k+1}$, therefore
$$(\eta+k(1-\theta))\log H\leq \log Z\leq (k+1)\log H\,.$$
Then $H^{1-\delta+1-\theta}\leq DY\leq H^{1+\delta+1}$, so
$$(2-( \delta+\theta))\log H\leq \log(DY)\leq (2+\delta)\log H\,.$$
Further from $H^\eta\leq AX\leq H^{1+k\delta+1}$ we get
$$\eta\log H\leq \log(AX)\leq (2+k\delta)\log H\,.$$
We finally get
$$\frac{\eta}{k+1}\log Z\leq\log(AX)\leq\frac{2+k\delta}{\eta+k(1-\theta)}\log Z$$
and 
\begin{equation}\label{eq:DYZ}\frac{2-(\delta+\theta)}{k+1}\log Z\leq\log(DY)\leq \frac{2+\eta}{\eta+k(1-\theta)}\log Z\,,
\end{equation}
so indeed condition (1) holds. Observe that we changed earlier the interval of definition for $x$ from $x\geq 1$ to $x\geq H^\eta$ namely because of this condition. In case $\eta=0$ we have $\log(AX)\geq 0$ and we cannot assure existence of a positive constant $C_1$ such that $C_1\log Z\leq \log(AX)$.\\

Condition (2) of Theorem \ref{thm:reuss} requires $l\geq 2$, or $DY\gg Z^{1/k}$. Our Diophantine equation is  $xy^k-ad^k=-C$ with $l=1$, so we need to verify $DY\gg Z^{1/k}$. Let us again impose a condition on the parameter $\delta$: $\delta+\theta<1-1/k$. Then from (\ref{eq:Det_int}) it follows  
$$DY\geq H^{2-(\delta+\theta)}>H^{1+1/k}\geq (XY^k)^{1/k}=Z^{1/k}$$
and in this case the second condition is also fulfilled. Therefore Theorem \ref{thm:reuss} does apply.\\

In such case we obtain
$$\mathcal{N}(Z;D,Y)\ll Z^\varepsilon\min\left\{(DYM)^{1/2}+D+Y, (AXM)^{1/2}+X+A\right\},$$
where $$\log M=\frac 9 8\frac{\log(DY)\log(AX)}{\log Z}\,.$$
Since $Z=XY^k=AD^k$, we have $Z^2=AX(DY)^k$ and $2\log Z=\log(AX)+k\log(DY)$. Then we can transform
$$\log M=\frac 9 8 \frac{\log(DY)}{\log Z}(2\log Z-k\log(DY))=\frac 9 8 \log(DY)\left(2-k\frac{\log(DY)}{\log Z}\right)\,.$$
From (\ref{eq:DYZ}) we see that 
$$\frac{2-(\delta+\theta)}{k+1}\leq\frac{\log(DY)}{\log Z}\,,$$
so we can bound from above
$$\log M\leq \frac 9 8 \log(DY)\left(2-k\frac{2-(\delta+\theta)}{k+1}\right)\leq  \frac 9 8 \left(2-k\frac{2-(\delta+\theta)}{k+1}\right)(2+\delta)\log H\,.$$
We denote 
\begin{equation}\label{def:g} g=g(k,\delta,\theta):=\frac 9 8 \left(2-k\frac{2-(\delta+\theta)}{k+1}\right)(2+\delta)\,.
\end{equation}
Then $M\leq H^g$ and 
\begin{eqnarray}\label{eq:N}\mathcal{N}(Z;D,Y) & \ll & H^{(k+1)\varepsilon}\left((DYM)^{1/2}+D+Y\right)\ll H^{(k+1)\varepsilon}\left(H^{1+\delta /2}H^{g/2}+H^{1+\delta}\right)\nonumber\\
&\ll & H^{\varepsilon}\left(H^{1+\delta/2+g/2}+H^{1+\delta}\right)\,.
\end{eqnarray}
%Let us choose $\varepsilon>0$ from Theorem \ref{thm:reuss} such that $(k+1)\varepsilon=\delta/2$. Then we get
%$$\mathcal{N}(Z;D,Y)\ll H^{1+\delta+g/2}+H^{1+3\delta/2}$$
Here $\varepsilon>0$ from Theorem \ref{thm:reuss} is arbitrary small and is present in the Vinogradov symbol dependence.\\ 

First we will assure that $1+\delta/2+g/2<2$ with a suitable choice of $\delta$ and $\theta>\delta$. Take 
$$\theta=2\delta\,.$$ 
Then 
\begin{eqnarray*}G(k,\delta)&:=&1+\frac \delta 2+\frac g 2 =  1+\frac \delta 2+\frac 9 8 \left(2-k\frac{2-3\delta}{k+1}\right)\left(1+\frac \delta 2\right)\\
%%&=&\left(1+\frac \delta 2\right)\left[1+\frac 9 8\left(2-\frac{2k}{k+1}+\frac{3k}{k+1}\delta\right)\right]\\
&=&\left(1+\frac \delta 2\right)\left[1+\frac 9 8\left(\frac{2}{k+1}+\frac{3k}{k+1}\delta\right)\right]
\end{eqnarray*}

From the form of the error terms in (\ref{eq:SOBL}) it is clear that we would be happy with as big $\delta$ as possible. Let us write $\delta=1/(\alpha k)$ for some integer $\alpha$. From the conditions on $\delta$ and $\theta$, which we want to hold for any $k\geq 2$, e.g. $\delta+\theta=3\delta<1-1/k$ and $\delta<1/(k+1)$, we should have $\alpha\geq 4$. We look at $G(k,1/(\alpha k))$.

\begin{eqnarray*}G\left(k,\frac 1 {\alpha k}\right) & = & \left(1+\frac 1 {2\alpha k}\right)\left[1+\frac 9 8\left(\frac{2}{k+1}+\frac{3k}{k+1}\frac 1 {\alpha k}\right)\right]\\
&=&\frac{2\alpha k+1}{2\alpha k}\frac 1 {8\alpha(k+1)}\left(8\alpha(k+1)+18\alpha+27\right)
\end{eqnarray*}

We want to choose $\alpha$ so that $2>G(k,1/(\alpha k))$, i.e. $2 (2\alpha k)(8\alpha(k+1))>(2\alpha k+1)(8\alpha k+26\alpha+27)$. After performing the operations in the expression we get 
$$32\alpha^2 k^2+32\alpha^2 k>16\alpha^2 k^2+52\alpha^2 k+62\alpha k+26\alpha+27\,.$$
We regroup the summands and arrive at the problem of finding parameter $\alpha$ for which the following quadratic form is positive for any $k\geq 2$:
\begin{equation}\mathcal{Q}(k,\alpha):=16\alpha^2 k^2-(20\alpha^2+62\alpha)k-(26\alpha+27)\,.
\end{equation}
We calculate the zeros of $\mathcal{Q}(k,\alpha)$ for a fixed $\alpha\in\left\{4,5,6,7\right\}$ and we see that the first value for which the larger zero of the quadratic equation is smaller than $2$ is $ \alpha=7$. That is why we choose $\delta=1/(7k)$, in which case $\mathcal{Q}(k,7)>0$ for any $k\geq 2$, i.e. $G(k,1/(7k))<2$. \\

More precisely we consider
\begin{equation}\label{eq:G}G_k:=G\left(k,\frac 1 {7k}\right)=\left(1+\frac 1 {14k}\right)\left(1+\frac{9\cdot 17}{7\cdot 8}\frac 1 {k+1}\right)\,.
\end{equation}
Then from (\ref{eq:dyadic}) and (\ref{eq:N}), and $1+1/(7k)<G_k$, it follows that
\begin{equation}\label{eq:S3final} \widetilde{S}_3=\Oc\left(H^{\varepsilon+G_k}\right)\,.
\end{equation}
Now in  (\ref{eq:SOBL}) we choose $\eta=1/5$. From (\ref{eq:SOBL}) and (\ref{eq:S3final}) we conclude that 
$$S(H)=c_{f}H^2+\Oc\left(H^{2-1/(7k)}\right)+\Oc\left(H^{2-1/14+1/(14k)}\right)+\Oc\left(H^{\varepsilon+G_k}\right)\,.$$ 
Note that for $k\geq 3$ we can write 
$$S(H)=c_{f}H^2+\Oc\left(H^{2-1/(7k)}\right)$$
because the dominating error term is the latter one. For $k=2$ we have $G_2>2-1/14+1/(14\cdot 2)>2-1/(7\cdot 2)$, so we content ourselves to write down
$$S(H)=c_{f}H^2+\Oc\left(H^{2-\delta}\right),$$
for some real $\delta=\delta(k,f)>0$. This proves Theorem \ref{thm:L1}.

%%----------------------------------------------------------------------

\section{Proof of Theorem \ref{thm:L2}}\label{sec:5}
The method used to attack Erd\H{o}s' conjecture from \cite{brown}, \cite{heath-br2} and \cite{reuss2} can be extended to the analogous conjecture for our two-variable polynomial. In a similar way like in (\ref{eq:S_split}) we split the sum $S'(H)$ into few parts. Define \begin{equation*}S'(m,H):=\sum_{\substack{p,q\leq H\\ f(p,q)\equiv 0(m)}}1\,.
\end{equation*}
We can write 
\begin{equation}\label{eq:S'_split} S'(H)=S'_1+S'_2+S'_3\,,
\end{equation}
where 
$$S'_1:=\sum_{d\leq w}\mu(d)S'(d^k,H)$$ 
for a positive parameter $w<H^{1/k}$ which we will choose very soon. Further 
$$S'_2:=\sum_{w<d\leq H^{1/k}}\mu(d)S'(d^k,H)\,.$$ 
We can directly estimate the last sum $S'_3$ using (\ref{eq:S2final}), (\ref{eq:S4final}), (\ref{eq:S31}), (\ref{eq:S32}) and (\ref{eq:S3final}):
\begin{equation}\label{eq:S'3final} S'_3:=\sum_{H^{1/k}<d\ll H^{1+1/k}}\mu(d)S'(d^k,H)\ll \sum_{H^{1/k}<d}S(d^k,H)\ll H^{2-\delta}\,,
\end{equation}
where $\delta>0$ is the constant from Theorem \ref{thm:L1}.

%%----------------------------------------------------------------------

\subsection{Estimation of $S'_1$}
We will derive the main term in the desired asymptotic formula from $S'_1$. Recall the definition of the function $\rho(m)$ from Theorem \ref{thm:L1} and let $(\mu_1,\nu_1),\ldots,(\mu_r,\nu_r)$ be the solutions of the congruence $f(\mu,\nu)\equiv 0\pmod{d^k}$, where $r=\rho(d^k)$. Then
\begin{eqnarray*} S'(d^k,H)&=&\sum_{i\leq r}\#\left\{(p,q)\in\left[1,H\right]^2:\, p\equiv \mu_i(d^k), \, q\equiv\nu_i(d^k)\right\}\\
&=&\sum_{i\leq r}\#\left\{p\leq H: p\equiv \mu_i(d^k)\right\}\cdot\#\left\{q\leq H: q\equiv\nu_i(d^k)\right\}\\
%&=&\sum_{i\leq r}\sum_{\substack{p\leq H\\p\equiv \mu_i(d^k)}}\sum_{\substack{q\leq H\\q\equiv \nu_i(d^k)}}1\\
&=&\sum_{\substack{i\leq r\\(\mu_i,d)=1\\(\nu_i,d)=1}}\pi(H;d^k,\mu_i)\pi(H;d^k,\nu_i)+\Oc\left(\rho(d^k)\right)\\
&+&\sum_{\substack{i\leq r\\(\mu_i,d)=1\\(\nu_i,d)>1}}\pi(H;d^k,\mu_i)\cdot\Oc(1)+\sum_{\substack{i\leq r\\(\nu_i,d)=1\\(\mu_i,d)>1}}\pi(H;d^k,\nu_i)\cdot\Oc(1)\,.
\end{eqnarray*}

If $d\leq (\log H)^N$ for arbitrary large $N$, then by Siegel-Walfisz theorem there exists a positive constant $c_N$, depending only on $N$, such that 
\begin{eqnarray*}S'(d^k,H)&=&\sum_{\substack{i\leq r\\(\mu_i,d)=1\\(\nu_i,d)=1}}\left[\frac{\pi(H)}{\varphi(d^k)}+\Oc\left(He^{-c_N\sqrt{\log H}}\right)\right]^2\\
&+&\Oc\left(\rho(d^k)\right)+\Oc\left(\rho(d^k)\left(\frac{\pi(H)}{\varphi(d^k)}+\Oc\left(He^{-c_N\sqrt{\log H}}\right)\right)\right)\,.
\end{eqnarray*}

By Lemma $1$(Hooley, \cite{hooley09}) we get $\rho'(d^k)\leq\rho(d^k)\ll d^{2k-2+\epsilon}$ for $k\geq 2$, and $\varphi(n)\gg_\varepsilon n^{1-\varepsilon}$ (Theorem $327$, \cite{hardy-wright}). So there is a positive constant $c>0$, which might vary in its different occurences, such that for $d\leq (\log H)^N$ we obtain
\begin{eqnarray}
\label{eq:S-W}S'(d^k,H)=\rho'(d^k)\frac{\pi(H)^2}{\varphi(d^k)^2}+\Oc\left(H^2e^{-c\sqrt{\log H}}\right).\end{eqnarray}
Let us then choose $w=(\log H)^N$ so that we apply (\ref{eq:S-W}). We get
\begin{eqnarray}\label{eq:S'err} S'_1&=&\sum_{d\leq w}\mu(d)S'(d^k,H)=\pi(H)^2\sum_{d\leq w}\frac{\mu(d)\rho'(d^k)}{\varphi(d^k)^2}+\Oc\left(wH^2e^{-c\sqrt{\log H}}\right)\nonumber\\
&=&c'_{f}\pi(H)^2+\Oc\left(\pi(H)^2\sum_{d> w}\frac{\rho'(d^k)}{\varphi(d^k)^2}\right)+\Oc\left(H^2e^{-c\sqrt{\log H}}\right)\,.
\end{eqnarray}

First of all we pay attention to the constant 
$$c'_{f}=\sum_{d=1}^\infty \frac{\mu(d)\rho'(d^k)}{\varphi(d^k)^2}\,.$$
By $\rho'(p^k)<\rho(p^k)\ll p^{2k-2+\epsilon}$ for $k\geq 2$, and $\varphi(n)\gg_\varepsilon n^{1-\epsilon'}$, we get 
$$c'_{f}=\sum_{d=1}^\infty \frac{\mu(d)\rho'(d^k)}{\varphi(d^k)^2}\ll \sum_{d=1}^\infty \frac{d^{2k-2+\epsilon}}{d^{2k-\epsilon'}}=\sum_{d=1}^\infty d^{-2+\epsilon}<\infty\,.$$

It is well-known that $\rho'(m)$ is multiplicative, so we can write
$$c'_{f}=\prod_p \left(1-\frac{\rho'(p^k)}{\varphi(p^k)^2}\right).$$

We also note that $c'_{f}\neq 0$. We need $\rho'(p^k)<\varphi(p^k)^2$, but even $\rho'(p^k)\leq\varphi(p^k)$ holds. Indeed, as $\rho'(p^k)$ counts $(\mu,\nu)\in(\Z/p^k\Z)^2$ such that $(\mu,p)=(\nu,p)=1$, for a fixed such $\nu$, which can take $\varphi(p^k)$ values, we have only one solution of the congruence $\mu\equiv -C\bar{\nu}^k\pmod{p^k}$.\\

Let us go back to the formula (\ref{eq:S'err}). %For the last error term we have
%\begin{equation}\label{eq:Err4}\sum_{\substack{d\leq w\\\mu(d)\neq 0}}\rho(d^k)\ll \sum_{d\leq w}d^{2k-2}\ll w^{2k-1}=\left(\log H\right)^{N(2k-1)}\,.
%\end{equation}
For the first error term we first use that
$$\sum_{d> w}\frac{\rho'(d^k)}{\varphi(d^k)^2}\leq \sum_{d> w}\frac{\rho(d^k)}{\varphi(d^k)^2}\ll \sum_{d>w}\frac{d^{2k-2+\epsilon}}{d^{2k-\varepsilon}}=
\sum_{d>w}d^{-2+\varepsilon}\ll w^{-1+\varepsilon}=\left(\log H\right)^{N(\varepsilon-1)}\,.$$
Then the error term is itself bounded from above by
\begin{equation}\label{eq:Err1}  \frac{H^2}{(\log H)^2}\cdot \frac{1}{\left(\log H\right)^{N(1-\varepsilon)}} \ll \frac{H^2}{(\log H)^K}\,,
\end{equation}
where $K>2$ is any real number, as long as we choose $N$ and $\varepsilon$ appropriately after we have already fixed $K$. For example we can choose $N=2K$. \\

For the second error term in (\ref{eq:S'err}) we clearly have
\begin{equation}\label{eq:Err2}H^2e^{-c\sqrt{\log H}}\ll \frac {H^2}{(\log H)^K}\,.\end{equation}
%Then the error term is bounded from above by
%\begin{equation}\label{eq:Err2}
%\frac{H^2}{\exp\left(c_N/2\sqrt{\log H}\right)}\cdot (\log H)^{N(k-1+\varepsilon)}\ll \frac{H^2}{(\log H)^K}
%\end{equation}
%
%For the third error term in (\ref{eq:S'err}) we basically use the upper bound (\ref{eq:Err4}) for the sum $\sum_{d\leq w}\rho'(d^k)$, so
%\begin{equation}\label{eq:Err3}
%\frac{H^2}{\exp\left(2c_N\sqrt{\log H}\right)}\cdot \left(\log H\right)^{N(2k-1)}\ll \frac{H^2}{(\log H)^K}
%\end{equation}

From the formula (\ref{eq:S'err}) and the error terms estimates (\ref{eq:Err1}) and (\ref{eq:Err2}) we conclude that for any $K>2$ we have
\begin{equation}\label{eq:S'1final} S'_1=\sum_{d\leq (\log H)^{2K}}\mu(d)S'(d^k,H)=c'_{f}\pi(H)^2+\Oc\left(\frac{H^2}{(\log H)^K}\right)\,.
\end{equation}

%%----------------------------------------------------------------------

\subsection{Estimation of $S'_2$} It is enough to estimate the sum $S'_2$ trivially. Recall that  we have already done similarly in (\ref{eq:SmH}). Applying again Lemma 1(Hooley, \cite{hooley09}) for square-free $d$, and using that $1\ll H/d^k\ll H^2/d^{2k}$ if $d\leq H^{1/k}$, we get
$$S(d^k,H)=\rho(d^k)\left[\frac{H^2}{d^{2k}}+\Oc\left(\frac H {d^k}\right)+\Oc(1)\right]\ll \rho(d^k)\frac{H^2}{d^{2k}}\ll H^2 d^{-2+\epsilon}\,.$$
Then 
\begin{equation} \label{eq:S'2final}
S'_2\ll H^2\sum_{(\log H)^{2K}<d\leq H^{1/k}}d^{-2+\epsilon}\ll H^2\cdot(\log H)^{2K(-1+\epsilon)}\ll \frac{H^2}{(\log H)^K}\,.
\end{equation}

Finally Theorem \ref{thm:L2} follows from formula (\ref{eq:S'_split}) and the estimates (\ref{eq:S'3final}), (\ref{eq:S'1final}) and (\ref{eq:S'2final}).
\section{Some speculations}
It is plausible that a similar asymptotic formula as in Theorem \ref{thm:L1} holds for the number of the $(d-1)$-free values of any irreducible two-variable polynomial $f(x,y)\in\Z\left[x,y\right]$ of degree $d$, which does not have a fixed $(d-1)$-th power prime divisor. With the same straightforward approach we get to the problem of estimating the number of solutions of the equation $f(x,y)=ab^{d-1}$, which is hard even in the case $x^ly^k+C=ab^{d-1}$ for general $l+k=d$. It is not clear if the determinant method can be used in such more general setting. 
%The polynomial $f(x,y)=x^ly^k+C$ of degree $d=k+l$ and its $m$-free values, (or $m=k$ ?)

\subsection*{Acknowledgments.} We thank Imre Ruzsa for posing the question if the square-free values of certain discriminants have a positive density. His question motivated this work. We would like to thank Andr\'as Bir\'o for the helpful discussions.
 \par The author is partially supported by OTKA grant no. K104183.

\end{document}